\documentclass[12pt]{article}
\usepackage{amssymb,amsmath,epsfig,amscd,eucal}

\newtheorem{thm}{Theorem}[section]
\newtheorem{defn}[thm]{Definition}
\newtheorem{lem}[thm]{Lemma}
\newtheorem{prop}[thm]{Proposition}

\newtheorem{cor}[thm]{Corollary}
\newtheorem{exa}[thm]{Example}

\newtheorem{qu}[thm]{Question}
\newtheorem{letterthm}{Main Theorem}

\newenvironment{pf}{\par\medskip\noindent{\em Proof. }}{\hfill $\square$\par\medskip}
\newenvironment{pfof}[1]{\par\medskip\noindent{\em Proof of #1. }}{\hfill $\square$\par\medskip}

\newcommand{\F}{\mathbb{F}}

\newcommand{\Z}{\mathbb{Z}}

\newcommand{\curlyT}{\mathcal{T}}

\title{Residually free 3-manifolds}
\author{Henry Wilton}

\date{5th September 2008}

\begin{document}
\maketitle

\begin{abstract}
We classify those compact 3-manifolds with incompressible toral boundary whose fundamental groups are residually free. For example, if such a manifold $M$ is prime and orientable and the fundamental group of $M$ is non-trivial then $M \cong \Sigma\times S^1$, where $\Sigma$ is a surface.
\end{abstract}

The theory of \emph{residually free groups}, and in particular \emph{fully residually free groups}, has been the subject of intense study in recent years, for the most part because of the important role that finitely generated fully residually free groups (also known as \emph{limit groups}) play in the logic of free groups (see \cite{sela_diophantine_2001} \emph{et seq.}, and \cite{kharlampovich_irreducible_1998-1}, \cite{kharlampovich_irreducible_1998} \emph{et seq.}).  A group $G$ is residually free if for each non-trivial element $g$ there is a homomorphism $f$ from $G$ to a free group so that $f(g)$ is non-trivial; it is fully residually free if the same holds when the element $g$ is replaced by a finite set of non-trivial elements.  The simplest examples of residually free groups, aside from free groups and free abelian groups, are the fundamental groups of (most) compact surfaces.

The techniques of residual \emph{finiteness} have been applied successfully in low-dimensional topology, and in particular in the study of 3-manifolds.  Indeed, Haken 3-manifolds have residually finite fundamental group \cite{hempel_residual_1987}, and attempts to prove stronger results about the pro-finite topology on certain 3-manifold groups are the subject of much active research.  There is also interest in \emph{large} 3-manifolds---those whose fundamental groups virtually surject non-abelian free groups.  It seems natural to speculate that the theory of residually free groups might also be fruitful in the 3-manifold context. In this spirit, Leonid Potyagailo asked which 3-manifold groups are residually free \cite{sela_diophantine_}.

At least for closed 3-manifolds, the theory of residually free groups is too powerful.  Our main theorem asserts that there are very few interesting examples of prime, compact 3-manifolds with incompressible toral boundary whose fundamental groups are residually free.

\begin{letterthm}\label{t:Theorem A}
If $M$ is a prime, compact 3-manifold with incompressible toral boundary then $\pi_1(M)$ is residually free and non-trivial if and only if $M$ is one of the following.
\begin{enumerate}
\item A trivial circle bundle over an orientable surface.
\item A circle bundle with trivial monodromy over a non-orientable surface of Euler characteristic less than -1.
\item The non-trivial circle bundle with trivial monodromy over the projective plane.
\end{enumerate}
\end{letterthm}
Free products of residually free groups are not always themselves residually free.  In corollary \ref{All 3-manifolds} we remove the assumption of primeness and classify all compact residually free 3-manifolds with incompressible toral boundary.

The article is organized as follows.  In section \ref{Rfg section} we recall the theory of residually free and fully residually free groups.  Subsection \ref{Products subsection} pays particular attention to their respective behaviours under free and direct products.  In a nutshell, residually free groups are well behaved under direct products but not under free products, whereas fully residually free groups are well behaved under free products but not under direct products.  In subsection \ref{Abelian splittings subsection} we prove necessary conditions for residually free groups to split over abelian subgroups.  Subsection \ref{ss: Virtual retractions} contains an elementary proof that residually free groups virtually retract onto cyclic subgroups.

Section \ref{3m section} is concerned with applying this theory in the context of compact 3-manifolds with incompressible toral boundary.  In subsection \ref{ss: Circle bundles}, we classify all residually free circle bundles. We then examine the torus decomposition of an arbitrary (compact, irreducible, orientable) 3-manifold with incompressible toral boundary.  We examine each piece and prove that, if residually free, it must be a circle bundle over a surface.  The conditions from subsection \ref{Abelian splittings subsection} are enough to deduce that any residually free graph manifold has trivial torus decomposition.

In section \ref{s: General boundary} we discuss the case of manifolds with more general boundary.  A complete classification here seems intractable without answers to some fundamental questions in the theory of fully residually free groups.

\subsection*{Acknowledgements}

Thanks to Cameron Gordon, Ben Klaff, Alan Reid and Matthew Stover for patient conversations, and to Martin Bridson for pointing out the alternative argument described in subsection \ref{ss: Bridson}.  Thanks also to the referee for pointing out several problems in an earlier draft.

\section{Residually free groups}\label{Rfg section}

\subsection{Free products and direct products}\label{Products subsection}

We shall start by recalling some of the theory of residually free and fully residually free groups.  Several of the results of this section are old results of \cite{baumslag_residually_1967}.  The point of view that we emphasize is that residually free groups are well behaved under direct products, whereas fully residually free groups are well behaved under free products.  Fix a non-abelian free group $\F$ of finite rank at least 2.

\begin{defn}
A group $G$ is \emph{residually free} if, for every non-trivial $g\in G$, there exists a homomorphism $f:G\to\F$ such that $f(g)\neq 1$.
\end{defn}
It is immediate that subgroups of residually free groups are residually free, and that direct products of free groups are residually free.  It is also clear that residually free groups are torsion-free.

\begin{exa}
Any subgroup of a direct product of free groups is residually free.
\end{exa}
We also have the following classification result.

\begin{lem}\label{2-generator rf groups}
Every 2-generator residually free group is isomorphic to one of $1$, $\Z$, $\Z^2$ or $F_2$, the non-abelian free group of rank 2.
\end{lem}
\begin{pf}
Let $G$ be a non-abelian residually free group generated by $a$ and $b$.  Since $[a,b]\neq 1$ there exists a homomorphism $f:G\to\F$ such that $f([a,b])\neq 1$.  The image of $f$ is a 2-generator, non-abelian subgroup of $\F$ and so is isomorphic to $F_2$.  Since $G$ is a 2-generator group that surjects onto $F_2$ it follows that $G\cong F_2$.
\end{pf}

We will need one consequence of this.  Let $G$ be a group and $H$ a subgroup. Denote by $Z_G(H)$ the centralizer of $H$, and $Z(G)=Z_G(G)$ the centre of $G$.  Let $[G,G]$ be the commutator subgroup of $G$.

\begin{lem}\label{Centralizer of the commutator subgroup}
If $G$ is residually free then
$$
Z_G([G,G])=Z(G).
$$
\end{lem}
\begin{pf}
Let $z\in Z_G([G,G])$.  For every $g\in G$, $[z,[z,g]]=1$. Therefore $\langle z,g\rangle$ is not isomorphic to $F_2$, and so $[z,g]=1$ by lemma \ref{2-generator rf groups}. Hence $Z_G([G,G])\subset Z(G)$.  Since the reverse inclusion is obvious, the result follows.
\end{pf}

In order to understand residually free groups, we need to consider the special class of \emph{fully residually free groups}.

\begin{defn}
A group $G$ is \emph{fully residually free} (also \emph{$\omega$-residually free} or, in the finitely generated case, a \emph{limit group}) if, for any finite collection of non-trivial elements $g_1,\ldots, g_n\in G$, there exists a homomorphism $f:G\to\F$ such that $f(g_i)\neq 1$ for all $1\leq i\leq n$.
\end{defn}

Of course, free groups are fully residually free. It is an easy exercise to show that free abelian groups are fully residually free.  The next examples of fully residually free groups are (most) surface groups.

\begin{defn}
If $\Sigma$ is a compact (not necessarily closed) surface and $\pi_1(\Sigma)$ is residually free then we call $\Sigma$ a \emph{residually free surface}.
\end{defn}
If $\Sigma$ has non-empty boundary then $\pi_1(\Sigma)$ is free, and so certainly fully residually free.

\begin{exa}
Let $\Sigma$ be a closed surface.  If $\Sigma$ is one of the non-orientable surfaces of Euler characteristic -1, 0 or +1 then $\pi_1(\Sigma)$ is \emph{not} residually free.  Otherwise, $\pi_1(\Sigma)$ is fully residually free.
\end{exa}

It is clear that the fundamental groups of the projective plane and the Klein bottle are not residually free.  Lyndon \cite{lyndon_equation_1959} proved that the fundamental group of the surface of Euler characteristic -1 is not residually free.  See lemma \ref{l: Folklore lemma} for an idea of how to start proving that the remaining surface groups are fully residually free.

The class of residually free groups is closed under taking direct products.  But a free product is rarely residually free.

\begin{lem}[Cf. theorem 6 of \cite{baumslag_residually_1967}]\label{Free products of rf groups}
If $A$ and $B$ are groups, $A$ is non-trivial and $B$ is not fully residually free then $A*B$ is not residually free.
\end{lem}
We prove this in the case when $B$ is not \emph{2-residually free}, that is when there exist non-trivial $b_1, b_2\in B$ such that every homomorphism $B\to\F$ kills one of $b_1$ or $b_2$.  The general proof is similar.  (In fact, every 2-residually free group is fully residually free, so this is the general case \cite{remeslennikov_exists-free_1989}.)

\begin{pfof}{lemma \ref{Free products of rf groups}}
Let $a\in A$ be non-trivial and let $b_1,b_2\in B\smallsetminus 1$ be such that every homomorphism $B\to\F$ kills one of $b_1$ or $b_2$. Define
$$
\xi=[b_1,[b_2,a]].
$$
Since
$$
\xi=b_1b_2ab_2^{-1}a^{-1}b_1^{-1}ab_2a^{-1}b_2^{-1}
$$
is a reduced word in the free product decomposition of $G$, it follows that $\xi\neq 1$.  But $f(\xi)=1$ for every homomorphism $f:A*B\to\F$, so $G$ is not residually free.
\end{pfof}

By contrast, the class of \emph{fully} residually free groups is well known to be closed under taking free products, but a direct product is rarely fully residually free.

\begin{lem}\label{Commutative transitivity}
If $G$ is fully residually free then $G$ is \emph{commutative transitive}; that is, if $a,b,c\in G\smallsetminus 1$ with $[a,b]=[b,c]=1$ then $[a,c]=1$.
\end{lem}

Two immediate consequences are that a non-abelian fully residually free group has trivial centre, and that if $A$ is non-trivial and $B$ is non-abelian then $A\times B$ is not fully residually free.  The proof of lemma \ref{Commutative transitivity} is well known; however, it will be important, so we recall it here.

\begin{pfof}{lemma \ref{Commutative transitivity}}
Suppose $[a,c]\neq 1$.  Then there exists a homomorphism $f:G\to\F$ so that $f(b)\neq 1$ and $f([a,c])\neq 1$.  But
$$
[f(a),f(b)]=[f(b),f(c)]=1
$$
so $f(a)$, $f(b)$ and $f(c)$ all lie in the same cyclic subgroup of $\F$; so $f([a,c])=1$, a contradiction.
\end{pfof}

In short: residually free groups are well behaved under direct products, but not under free products; by contrast, fully residually free groups are well behaved under free products, but not direct products.

From this point of view, the simplest residually free group that is not fully residually free is $\F\times\Z$.  Indeed, it is a theorem of B.~Baumslag that this is the only obstruction to being fully residually free.

\begin{thm}[Theorem 1 of \cite{baumslag_residually_1967}]\label{Baumslag's theorem}
Every finitely generated residually free group is either fully residually free or contains $\F\times\Z$.
\end{thm}

\subsection{Amalgamated products and HNN-extensions}\label{Abelian splittings subsection}

We have seen that fully residually free groups can be non-trivial free products, whereas residually free groups that are not fully residually free cannot be.  We now turn our attention to splittings (as amalgamated products or HNN-extensions) over abelian subgroups.

Theorem 3.2 of \cite{sela_diophantine_2001} assert that every finitely generated fully residually free group splits non-trivially over an abelian (possibly trivial) subgroup.  The following is a straightforward consequence.

\begin{thm}[\cite{sela_diophantine_2001}]\label{Splittings of frf groups}
Every finitely generated non-cyclic fully residually free group has infinitely many outer automorphisms.
\end{thm}

We shall now address splittings of residually free groups over abelian subgroups.  The idea is to generalize the proof of lemma \ref{Free products of rf groups}.  For any subgroup $H\subset G$, define the \emph{normal core} of $H$ in $G$ to be the largest normal subgroup of $G$ contained in $H$; that is,
$$
\mathrm{Core}_G(H)=\bigcap_{g\in G} H^g.
$$
Of course, $\mathrm{Core}_G(H)=H$ if and only if $H$ is normal in $G$.

\begin{prop}\label{Amalgamated product}
Let $G=A*_C B$ be an amalgamated product over an abelian subgroup, and suppose that $G$ is residually free. If $A$ is non-abelian then
$$
Z(A)\cap C\subset \mathrm{Core}_B(C).
$$
\end{prop}
\begin{pf}
Fix $z\in Z(A)\cap C$. Since $A$ is non-abelian, $A$ contains a non-abelian free subgroup by lemma \ref{2-generator rf groups}.  The edge group $C$ intersects this free subgroup in a cyclic subgroup, and so there exist $a_1,a_2\in A$ such that $[a_1,a_2]$ is not contained in $C$.  Let $b\in B\smallsetminus C$. Consider
$$
\xi = [[a_1,a_2],[z,b]].
$$
For every $f:G\to\F$, either  $f(z)=1$ or $f([a_1,a_2])=1$, (to see this, apply the proof of lemma \ref{Commutative transitivity} to the triple $a_1,z,a_2$) and so $f(\xi)=1$. Since $G$ is residually free, $\xi=1$.  Therefore
$$
\xi=[a_1,a_2][z,b][a_2,a_1][b,z]
$$
is not a reduced word in the amalgamated product decomposition of $G$, so $[z,b]\in C$.  Therefore every $b\in B$ conjugates $z$ into $C$. Hence $z\in\mathrm{Core}_B(C)$.
\end{pf}

In the case of an HNN-extension we get a stronger result, that holds for \emph{any} edge group.

\begin{prop}\label{HNN-extension}
Let $G=A*_C$ be an HNN-extension and suppose $G$ is residually free.  Let $t$ be a stable letter. Then
$$
(Z(A)\cap C)^t\subset Z(A).
$$
\end{prop}
\begin{pf}
Fix $z\in Z(A)\cap C$, and let $a_1,a_2\in A$. Consider
$$
\xi = [[a_1,a_2],[z,t]].
$$
For every $f:G\to\F$, either  $f(z)=1$ or $f([a_1,a_2])=1$ (again by the proof of lemma \ref{Commutative transitivity}) so $f(\xi)=1$. Since $G$ is residually free, $\xi=1$.  If $z^t=a$ then $\xi$ can be written as
$$
\xi=[[a_1,a_2],za^{-1}]=[[a_1,a_2],a^{-1}].
$$
So $a$ commutes with every commutator in $A$; that is, $a\in Z_A([A,A])$ which equals $Z(A)$ by lemma \ref{Centralizer of the commutator subgroup}.
\end{pf}

\subsection{Virtual retractions onto cyclic subgroups}\label{ss: Virtual retractions}

A subgroup $H$ of a group $G$ is a \emph{retract} if the inclusion map $H\hookrightarrow G$ has a left inverse.  Likewise $H$ is a \emph{virtual retract} of $G$ if it is contained as a retract in a finite-index subgroup of $G$.  In \cite{bridson_subgroup_2008} it is shown that residually free groups virtually retract onto many subgroups.  Here we give an elementary argument to prove the special case that cyclic subgroups of residually free groups are virtual retracts.  We shall apply this in subsection \ref{ss: Torus bundles} to prove that residually free torus bundles are 3-tori.

\begin{lem}\label{LR over Z}
Let $G$ be a residually free group and $Z$ a cyclic subgroup.  Then $G$ has a finite-index subgroup $G'$, containing $Z$, such that $Z$ is a retract of $G'$.
\end{lem}
\begin{pf}
Let $f:G\to\F$ be a homomorphism that is injective on $Z$.  By Marshall Hall's Theorem \cite{hall_subgroups_1949} there exists a finite-index subgroup $F'$ of $\F$ that retracts onto $f(Z)$.  Let $G'=f^{-1}(F')$.  Then $f$ composed with the retraction $F'\to f(Z)$ defines a retraction $G'\to Z$, as required.
\end{pf}

\section{3-manifolds}\label{3m section}

\subsection{Circle bundles}\label{ss: Circle bundles}

The observations of section \ref{Rfg section} give us some easy examples of 3-manifolds with residually free fundamental groups.

\begin{exa}
If $\Sigma$ is a residually free surface then $\Sigma\times S^1$ has residually free fundamental group.
\end{exa}

In this subsection we shall examine the question of which other circle bundles over closed surfaces have residually free fundamental group.  Any circle bundle $M$ over a surface $\Sigma$ is determined by its monodromy and by the Euler class.  If the fundamental group of $M$ is infinite then it fits into a short exact sequence
\[
1\to \Z\to\pi_1(M)\stackrel{\eta}{\to}\pi_1(\Sigma)\to 1.
\]
The action of $\pi_1(\Sigma)$ on $\Z$ by conjugation is the \emph{monodromy} of the bundle.  The Euler class is an element of $H^2(\Sigma,\Z)$, where the coefficients are twisted by the monodromy, and the Euler class itself can be interpreted as the obstruction to the sequence splitting. (See pages 434-5 of \cite{scott_geometries_1983} for more details.)

If the fundamental group of a circle bundle is residually free then it follows from lemma \ref{2-generator rf groups} that the monodromy is trivial.

\begin{exa}[Bundles over orientable surfaces]
Suppose $\Sigma$ is closed and orientable and $\pi_1(M)$ is non-trivial and residually free.  Then the fundamental group of $M$ has a presentation of the form
\[
\langle a_1,\ldots,a_g,b_1,\ldots,b_g,z\mid \prod_i[a_i,b_i]=z^e,[a_1,z]=\ldots=[b_g,z]=1\rangle
\]
where $e$ times the fundamental class of $\Sigma$ can be identified with the Euler class (up to sign).  Consider a homomorphism $f:\pi_1(M)\to\F$. If $f(z)$ is non-trivial then it commutes with everything else in the image, so the image is abelian and $f$ factors through $H_1(M,\Z)$.   If $e\neq 0$ then the homology class of $z$ is torsion, so $z$ dies in every map to $\F$.  Therefore only trivial circle bundles over orientable surfaces are residually free.

Because every other case has non-trivial centre but is not abelian, $\pi_1(M)$ is fully residually free only if $\Sigma$ is the 2-sphere or the 2-torus.
\end{exa}

We now consider the case in which $\Sigma$ is non-orientable.

\begin{exa}[Bundles over non-orientable surfaces]
Suppose $\Sigma$ is closed and non-orientable and $\pi_1(M)$ is residually free.  In this case $H^2(\Sigma,\Z)$ is isomorphic to $\Z/2$, so there are only two choices for the Euler class, and $M$ has a presentation of the form
\[
\langle a_1,\ldots,a_m,z\mid \prod_i a_i^2=z^e,[a_1,z]=\ldots=[b_m,z]=1\rangle
\]
where $e$ can be taken to be $0$ or $1$, depending on the Euler class.  In either case, it is clear that the infinite cyclic subgroup generated by $z$, which corresponds to the circle fibre, injects under the map to $H_1(M,\Z)$.

As $\langle z\rangle$ is precisely the kernel of the map $\eta:\pi_1(M)\to\pi_1(\Sigma)$, it follows that if $\Sigma$ is residually free (that is, if $\chi(\Sigma)<-1$) then $M$ has residually free fundamental group.

If $\Sigma$ has Euler characteristic $1$, $0$ or $-1$ then $\pi_1(\Sigma)$ is not residually free.  Indeed, in these cases it follows from \cite{lyndon_equation_1959} that every homomorphism $\pi_1(\Sigma)\to\F$ factors through $H_1(\Sigma,\Z)$.\footnote{It is an entertaining exercise to prove this topologically.  Consider a cellular map from $\Sigma$ to a graph whose fundamental group if $\F$, and pull back the midpoints of the edges. The resulting simple closed curves cut $\Sigma$ into simpler pieces.}  Just as in the orientable case, every homomorphism $\pi_1(M)\to\F$ either factors through $H_1(M,\Z)$ or $\pi_1(\Sigma)$, so in fact every homomorphism to $\F$ has abelian image.  If $\chi(\Sigma)$ is $0$ or $-1$, in which case $\pi_1(M)$ is not abelian, it follows that $\pi_1(M)$ is not residually free.  However, in the case when $\Sigma$ is the projective plane and $e=1$, we see that $\pi_1(M)$ is infinite cyclic, so residually free.

This last example is the only one that is fully residually free---all the others have a non-trivial centre but are not abelian.
\end{exa}

We summarize these examples in the following proposition.

\begin{prop}\label{p: Circle bundles}
If $M$ is a circle bundle with incompressible boundary over a compact surface $\Sigma$ then $\pi_1(M)$ is residually free and non-trivial if and only if $M$ is of one of the following forms.
\begin{enumerate}
\item The base $\Sigma$ is orientable and residually free and $M$ is the trivial bundle.
\item The base $\Sigma$ is non-orientable and residually free and $M$ is any circle bundle with trivial monodromy.
\item The base $\Sigma$ is the projective plane and $M$ is the non-trivial bundle with trivial monodromy.
\end{enumerate}
Of these, only $S^2\times S^1$, the 3-torus and the non-trivial bundle over the projective plane have fully residually free fundamental groups.
\end{prop}

Note that the third example can also be thought of as the twisted sphere bundle over the circle.

\subsection{The torus decomposition}

We shall apply the theory of residually free groups to a compact, irreducible, orientable 3-manifold $M$ with (possibly empty) toral boundary.  Our principal tool is the torus decomposition of Jaco--Shalen \cite{jaco_seifert_1979} and Johannson \cite{johannson_homotopy_1979}.

\begin{thm}\label{JSJ decomposition}
Let $M$ be a compact, irreducible, orientable 3-manifold with incompressible boundary.  Then $M$ can be cut along a canonical embedded collection of disjoint incompressible 2-sided tori
$$
\curlyT=T_1\sqcup\ldots\sqcup T_k
$$
such that each component of the complement $M\smallsetminus \curlyT$ is either atoroidal or Seifert-fibred.
\end{thm}

See \cite{neumann_canonical_1997} for an excellent exposition of this result. For (sketchy) definitions of atoroidal and Seifert-fibred manifolds, see the next two subsections.

Since being residually free passes to subgroups, we start by classifying the possible pieces of the torus decomposition of $M$.  If $\pi_1(M)$ is residually free and non-trivial then it clearly has positive first Betti number.  It follows that $M$ is Haken and so satisfies Thurston's geometrization conjecture, which asserts that each piece of the torus decomposition carries one of eight model geometries.  See \cite{scott_geometries_1983} for a detailed account.  In particular, we assume that the atoroidal pieces are either Seifert-fibred or carry hyperbolic structures.  Under the additional hypothesis that the boundary of $M$ is toral, these hyperbolic structures are of finite volume.  Although we give a unified treatment for the Seifert-fibred pieces in subsection \ref{Seifert-fibred pieces subsection}, the reader may find it instructive to consider which of the fibred geometries can be supported by a residually free 3-manifold.

\subsection{Hyperbolic pieces}\label{ss: Atoroidal pieces}

First, consider an \emph{atoroidal} piece.  This means that every essential embedded torus is boundary parallel.  As mentioned above, we assume that such pieces are Seifert-fibred or admit hyperbolic geometries of finite volume.  Seifert-fibred pieces are dealt with in the next subsection.  The following proposition rules out hyperbolic pieces.

\begin{prop}\label{Atoroidal pieces}
If $M$ is a hyperbolic 3-manifold of finite volume then $\pi_1(M)$ is not residually free.
\end{prop}
\begin{pf}
Suppose $\pi_1(M)$ is residually free.  By theorem \ref{Baumslag's theorem}, either $\pi_1(M)$ contains $\F\times\Z$ or $\pi_1(M)$ is fully residually free.  As $M$ is hyperbolic any $\Z^2$ subgroup of $\pi_1(M)$ is malnormal, and so $\pi_1(M)$ cannot contain $\F\times\Z$---therefore, it is fully residually free, and by theorem \ref{Splittings of frf groups} has infinitely many outer automorphisms.  But it is a well known consequence of Mostow Rigidity that the fundamental groups of finite-volume hyperbolic 3-manifolds have finite outer automorphism group.
\end{pf}

\subsection{Seifert-fibred pieces} \label{Seifert-fibred pieces subsection}

As none of the pieces of the torus decomposition can be hyperbolic, every compact 3-manifold with incompressible toral boundary and non-trivial, residually free fundamental group is a \emph{graph manifold}---that is, every piece of the torus decomposition is Seifert-fibred.  (For us, torus bundles over circles are graph manifolds.) Furthermore, the Seifert-fibred pieces must be of a very restricted type.

For details of the theory of Seifert-fibred manifolds see, for example, \cite{scott_geometries_1983}.  We will not define Seifert-fibred manifolds here, but we recall an elementary fact.  If $M$ is a Seifert-fibred manifold then $\pi_1(M)$ fits into a short exact sequence
$$
1\to Z\to\pi_1(M)\stackrel{\eta}{\to}\pi_1(\Sigma)\to 1
$$
where $Z$ is cyclic.  The following lemma reduces the Seifert-fibred pieces to the cases studied in subsection \ref{ss: Circle bundles}.

\begin{lem}\label{l: Seifert-fibred pieces}
Let $M$ be a Seifert-fibred manifold with incompressible boundary.  If $\pi_1(M)$ is non-trivial and residually free then $M$ has no singular fibres and so is a circle bundle over a surface.
\end{lem}
\begin{pf}
For a contradiction, assume that the base orbifold $\Sigma$ has non-trivial singular locus.  Because the fundamental group of the Klein bottle is not residually free, it follows that $\Sigma$ has no reflector lines and therefore the singular fibres of $M$ are isolated.  By lemma \ref{2-generator rf groups}, $Z$ is central in $\pi_1(M)$.   Let $g\in\pi_1(M)$ be an element corresponding to a singular fibre.  The image curve in $\Sigma$ can be taken to be the generator corresponding to a cone point, so $\eta(g)^p=1$ for some integer $p$.

It follows from the fact that $\pi_1(M)$ is infinite that the base orbifold $\Sigma$ has a manifold cover (see \cite{scott_geometries_1983}, lemma 3.1).  Because the boundary of $M$ is incompressible, $\Sigma$ is not a disc with a single cone point and hence $\eta(g)$ is not central in $\pi_1(\Sigma)$.  So there exists $h\in\pi_1(M)$ such that $\eta(g)$ and $\eta(h)$ do not commute.  But $[g^p,h]=1$, so $\langle g,h\rangle$ is abelian by lemma \ref{2-generator rf groups}.  This contradicts the fact that $\eta(g)$ and $\eta(h)$ do not commute.
\end{pf}

\subsection{Torus bundles}\label{ss: Torus bundles}

Torus bundles over circles are important special cases of graph manifolds.  If $M$ is such a torus bundle then its fundamental group fits into a short exact sequence
\[
1\to A\to\pi_1(M)\to\Z\to 1
\]
where $A\cong\Z^2$.  We will use lemma \ref{LR over Z} to show that any residually free torus bundle is actually Seifert-fibred.  The proof avoids the deep results of \cite{bridson_subgroup_2008}.

\begin{lem}\label{l: Torus bundles}
If $M$ is a torus bundle over the circle and $\pi_1(M)$ is residually free then $M$ is homeomorphic to a 3-torus.
\end{lem}
\begin{pf}
First, we will show that $M$ is Seifert-fibred.  It suffices to show that the action of $\Z$ on $A$ fixes some non-trivial element.  Consider an arbitrary element $a\in A$.  By lemma \ref{LR over Z}, there exists a finite-sheeted cover $M'$ of $M$ with a retraction $\rho:\pi_1(M')\to\langle a\rangle$.   Let $b$ generate $A\cap\pi_1(M')\cap\ker\rho$, so $A\cap\pi_1(M')=\langle a\rangle\times\langle b\rangle$.  For any $g\in\pi_1(M)$ there is an $m$ such that $g^m\in\pi_1(M')$.  Then
\[
\rho(g^m bg^{-m})=\rho(g^m)\rho(b)\rho(g^{-m})=1
\]
so $g^m bg^{-m}\in A\cap\pi_1(M')\cap\ker\rho$ and hence equals $b^n$ for some $n$.  Hence the subgroup of $\pi_1(M)$ generated by $g$ and $b$ is abelian by lemma \ref{2-generator rf groups}.   It follows that the action of $\Z$ on $A$ fixes $b$.

Therefore $M$ is Seifert-fibred.  So, by lemma \ref{l: Seifert-fibred pieces}, $M$ is a circle bundle over a surface and hence by proposition \ref{p: Circle bundles} a 3-torus.
\end{pf}

\subsection{Graph manifolds}

We have shown that every (closed, orientable, irreducible) residually free 3-manifold is a graph manifold, and that every Seifert-fibred piece is a circle bundle over a surface.

\begin{prop}\label{Graph manifolds}
Let $M$ be a compact, irreducible, orientable graph manifold with incompressible boundary.  If $\pi_1(M)$ is non-trivial and residually free then $M$ is Seifert-fibred.
\end{prop}
\begin{pf}
Consider an essential torus $T$ in the torus decomposition of $M$, and let $M'=M\smallsetminus T$.  There are two inclusions $T\hookrightarrow M'$ and each induces a foliation of $T$ by circles, inherited from the Seifert-fibred submanifolds of $M'$ that contain $T$.  We shall show that these two foliations coincide---hence the Seifert structures extend over $T$.  There are two cases.

Suppose $T$ is contained in disjoint Seifert-fibred submanifolds $M_1$ and $M_2$. Both $M_1$ and $M_2$ are products, by lemma \ref{l: Seifert-fibred pieces}: for $i=1$ or $2$, let $M_i=\Sigma_i\times S^1_i$, where $\Sigma_i$ is a surface with non-trivial boundary and $\chi(\Sigma_i)<0$, and $S^1_i$ is the circle factor of $M_i$. The centre of $\pi_1(M_i)$ is $\pi_1(S^1_i)$, and also
$$
\mathrm{Core}_{\pi_1(M_i)}(\pi_1(T))=\pi_1(S^1_i).
$$
It follows by proposition \ref{Amalgamated product} that the gluing along $T$ sends circle factors to circle factors.

Similarly, suppose each inclusion of $T$ is contained in the same Seifert-fibred submanifold $M'$.  Then  $M'=\Sigma'\times S^1$, where $\Sigma'$ is a surface with non-trivial boundary and non-positive Euler characteristic.  Let $t$ be a stable letter of the corresponding decomposition of $\pi_1(M)$ as an HNN-extension.  If $\Sigma'$ is an annulus then $M$ is a torus bundle over a circle and the result follows from lemma \ref{l: Torus bundles}.  Otherwise, $\chi(\Sigma')<0$ and the centre of $\pi_1(M')$ is $\pi_1(S^1)$, so $(\pi_1(S^1)\cap\pi_1(T))^t\subset\pi_1(S^1)$ by proposition \ref{HNN-extension}.  Therefore, the gluing along $T$ sends the circle factor of $M'$ to itself.

As the Seifert structure extends over $T$, it follows that $T$ was not really in the torus decomposition of $M$.  Hence $M$ has trivial torus decomposition, so is Seifert-fibred.
\end{pf}

\subsection{Conclusion}

Combining these results we achieve a complete classification of prime, compact 3-manifolds (with incompressible toral boundary) with residually free fundamental group.

\begin{thm}\label{Prime 3-manifolds}
If $M$ is a prime, compact 3-manifold with incompressible toral boundary and $\pi_1(M)$ is residually free and non-trivial then $M$ is a circle bundle over a surface, and hence of one of the forms listed in proposition \ref{p: Circle bundles}.
\end{thm}
\begin{pf}
Let $\hat{M}$ be a finite-sheeted orientable cover of $M$.  We first argue that $\hat{M}$ is Seifert-fibred.  If $\hat{M}$ is not irreducible then it is $S^2\times S^1$.  Otherwise $\hat{M}$ has a well-defined torus decomposition, and because $b_1(\hat{M})$ is positive we can take the pieces that are not Seifert-fibred to admit hyperbolic structures of finite volume.  By proposition \ref{Atoroidal pieces} there are no hyperbolic pieces, so $\hat{M}$ is a graph manifold and so, by proposition \ref{Graph manifolds}, is in fact Seifert-fibred.  Therefore $M$ is Seifert-fibred and so, by lemma \ref{l: Seifert-fibred pieces}, is a circle bundle over a surface.
\end{pf}

This completes the proof of the main theorem.  We saw above that, of all residually free groups, only the \emph{fully} residually free ones are well-behaved under taking free products.  It is therefore easy to see when connected sums are residually free.  Note that, for the following result, we need to assume the Poincar\'e Conjecture.

\begin{cor}\label{All 3-manifolds}
If $M$ is any compact 3-manifold with incompressible toral boundary and $\pi_1(M)$ is residually free and non-trivial then one of the following holds:
\begin{enumerate}
\item $\pi_1(M)$ is fully residually free and $M$ is a connected sum of finitely many copies of $S^2\times S^1$, the 3-torus, and the non-trivial circle bundle over the projective plane;
\item $\pi_1(M)$ is \emph{not} fully residually free and $M$ is one of the non-fully residually free circle bundles listed in proposition \ref{p: Circle bundles}.
\end{enumerate}
\end{cor}
\begin{pf}
Consider the Kneser--Milnor prime decomposition
$$
M=M_1\#\ldots\# M_n
$$
of $M$, where no $M_i$ is a 3-sphere.  The fundamental group of $M$ is a free product of the fundamental groups of the prime pieces $\{M_i\}$.  Assuming the Poincar\'e Conjecture, every prime piece has non-trivial fundamental group.

If one of the $M_i$ is not fully residually free then, by lemma \ref{Free products of rf groups}, $M=M_i$ and so we are in case 2.

Otherwise, $\pi_1(M_i)$ is fully residually free for each $i$.  By theorem \ref{Prime 3-manifolds} each $M_i$ is a circle bundle over a surface, and by proposition \ref{p: Circle bundles} each $M_i$ is of one of the three types listed.
\end{pf}

In summary, there are no interesting examples of closed 3-manifolds with residually free fundamental group.  What is the next best thing that one might hope for?  As far as residual properties go, the class of torsion-free nilpotent groups seems a natural intermediate class between free groups and finite groups.

\begin{qu}
Which (closed) 3-manifolds have fundamental groups that are residually torsion-free nilpotent?
\end{qu}

It is a theorem of Magnus (see \cite{lyndon_combinatorial_1977}, Proposition 10.2) that free groups are residually torsion-free nilpotent.  Furthermore, the classes of graph manifolds and right-angled Artin groups intersect in the class of \emph{right-angled tree groups} \cite{behrstock_quasi-isometric_2008}.  In \cite{duchamp_lower_1992} it is proved that all right-angled Artin groups are residually torsion-free nilpotent.  So the class of residually torsion-free nilpotent 3-manifolds is certainly substantially larger than the class of residually free 3-manifolds.

\subsection{An alternative argument}\label{ss: Bridson}

Martin Bridson has pointed out an alternative proof of the main theorem, which we outline here.  Let $M$ be an irreducible 3-manifold with non-trivial residually free fundamental group.  By the main theorem of \cite{bridson_subgroups_}, $\pi_1(M)$ is virtually a direct product of fully residually free groups.  If there were more than one non-abelian factor then, by lemma \ref{2-generator rf groups}, $\pi_1(M)$ would contain a subgroup isomorphic to $\F\times\F$, which is well known to contain subgroups that are finitely generated but not finitely presentable---this contradicts the fact that 3-manifolds are coherent \cite{scott_finitely_1973}, so at most one of the direct factors can be non-abelian.

There are now two cases to consider.  If there is a non-trivial abelian factor, then $\pi_1(M)$ has an infinite cyclic normal subgroup and so is Seifert-fibred; one then argues as in subsection \ref{Seifert-fibred pieces subsection}.  On the other hand, if there is no non-trivial abelian factor then the fundamental group is fully residually free and hence, by theorem \ref{Splittings of frf groups}, has infinite outer automorphism group. As in subsection \ref{ss: Atoroidal pieces}, this leads to a contradiction.  Bridson's argument is shorter than the one given here, but uses considerably more sophisticated machinery.

\section{More general boundary}\label{s: General boundary}

In this section we discuss compact 3-manifolds with arbitrary boundary.  We may assume that the boundary is still incompressible---otherwise the fundamental group is a free product, and so is easily dealt with by the results of section \ref{Rfg section}.

With arbitrary boundary, one can no longer guarantee that the atoroidal pieces of the torus decomposition admit hyperbolic structures of finite volume, so one can no longer apply theorem \ref{Splittings of frf groups} and Mostow Rigidity to rule them out.  There are certainly more examples of both residually free 3-manifolds in this context.  The following folklore lemma is useful for constructing fully residually free groups.

\begin{lem}\label{l: Folklore lemma}
Let $L$ be a fully residually free group and let $A\subset L$ be a maximal abelian subgroup.  The double $D=L*_A L$ is also fully residually free.
\end{lem}

Note that it follows from lemma \ref{Commutative transitivity} that if $A$ is not maximal in $L$ then $D$ is not fully residually free.

The idea of the proof of lemma \ref{l: Folklore lemma} is as follows.  Let $\rho: D\to L$ be the natural retraction that identifies the two copies of $L$ and let $\delta$ be an automorphism of $D$ given by a Dehn twist in some non-trivial element of $A$.  Then, for any element $g\in D$, $\rho\circ\delta^n(g)$ is non-trivial for all sufficiently large $n$.  (See \cite{bestvina_notes_} for some details.)

This enables us to construct examples of fully residually free 3-manifolds by gluing interval bundles along annuli.

\begin{exa}
Let $N$ be an interval bundle over a closed, (fully) residually free, hyperbolic surface $\Sigma$, and let $A$ be an essential annulus embedded in the boundary of $N$ whose core curve does does not have a proper root in $\pi_1(N)$.  Let $M$ be the double of $N$ along $A$.  Then $\pi_1(M)$ is fully residually free.
\end{exa}

Using lemma \ref{l: Folklore lemma} one can iterate the doubling construction to build still more complicated (fully) residually free 3-manifolds.  Of course, there are also negative examples, such as an interval bundle over the surface of Euler characteristic -1.  It seems hard to classify residually free manifolds in this context without a complete solution to the following challenge of Potyagailo and Sela.

\begin{qu}[\cite{sela_diophantine_}]
Give sufficient or necessary conditions for a group of the form $\F*_{u=v}\F$ to be fully residually free.
\end{qu}

Note that, by theorem \ref{Baumslag's theorem}, such a group is residually free if and only if it is fully residually free.

\bibliographystyle{plain}

\bigskip\bigskip\centerline{\textbf{Author's address}}

\smallskip\begin{center}\begin{tabular}{l}%
Department of Mathematics\\
1 University Station C1200\\
Austin, TX 78712-0257\\
USA\\
{\texttt{henry.wilton@math.utexas.edu}}
\end{tabular}\end{center}

\end{document}